\documentclass[12pt,dvips]{amsart}
\usepackage{pstricks, pstricks-add, pst-node}
\usepackage{amsmath,amstext,amsthm,amsfonts,latexsym,color}
\usepackage{graphicx}
\usepackage{graphics}
\newtheorem{teo}{Theorem}[section]

\newtheorem{lema}[teo]{Lemma}
\newtheorem{cor}[teo]{Corollary}
\newtheorem{maintheorem}{Theorem}
\newtheorem{mainlemma}{Lemma}

\newcommand{\bmt}{\begin{maintheorem}}
\newcommand{\emt}{\end{maintheorem}}
\newcommand{\bml}{\begin{mainlemma}}
\newcommand{\eml}{\end{mainlemma}}
\numberwithin{equation}{section}

\newcommand{\dist}{{\rm dist}}

\newcommand{\p}{\varphi}

\newcommand{\la}{\lambda}

\newcommand{\al}{\alpha}

\newcommand{\si}{\sigma}
\newcommand{\de}{\delta}

\newcommand{\R}{\mathbb{R}}
\newcommand{\Z}{\mathbb{Z}}
\newcommand{\T}{\mathbb{T}}

\newcommand{\m}{{\rm Leb}\:}

\newcommand{\fu}{\mathcal{F}^u}

\newcommand{\pern}{{\rm Per}_n}
\newcommand{\perti}{{\tilde{\rm Per}}_n}

\newcommand{\diam}{{\rm diam}\:}

\title[Bowen measure]{Bowen measure for derived from Anosov diffeomorphims}\label{SAMVA.tex}

\date{\today}
\author{Martin Sambarino}
\address{Martin Sambarino,
CMAT-Facultad de Ciencias,
Montevideo-Uruguay.
}
\email{samba@cmat.edu.uy}

\author{Carlos H. V\'asquez}
\address{Carlos H. V\'asquez,
Instituto de Matem\'atica,
Pontificia Universidad Cat\'olica de Valpara\'iso,
Valpara\'iso-Chile
}
\email{carlos.vasquez@ucv.cl}
\begin{thanks}{C.V. was partially supported by Fondecyt 11060497 and Research Network on Low Dimensional
Dynamics, PBCT ACT-17, CONICYT, Chile.}
\end{thanks}



\begin{document}

\maketitle

\begin{abstract}
In this work we give general conditions under which a $C^0$ perturbation of an expansive homeomorphim with specification property has an unique Bowen measure, that is, there is an ergodic probability measure which is the unique measure maximizing the topological entropy.

We apply these conditions to show that several derived from Anosov diffeomorphims have a unique Bowen measure.
\end{abstract}

\section{Introduction}\label{sec:intro}
The well known variational principle states that if $f:X\to X$ is a continuous map of the compact metric space $X$ then the supreme of the metric entropies over the set of invariant probabilities measure is the topological entropy:
$$h_{top}(f)=\sup_{\mu\in M(X,f)}h_\mu(f)$$
where $M(X,f)$ is the set of invariant probability measures.
It is an old and a very important problem  to know whether this supreme is in fact a maximum and in this case whether the maximum is attached by a single probability measure.

After Bowen, some sufficient conditions are known to give a positive answer to this question. Indeed he showed a way to construct it.
Let $X$ be a compact metric space and $f\::X\to X$ be an $\al$-expansive homeomorphism with the specification property. For $n\geq 1$, denote by
$$\pern(f):=\{x\in X\::\: f^ n(x)=x\},$$
and define the sequence of $f$-invariants measures
$$\mu_n:=\frac{1}{|\pern(f)|}\sum_{x\in \pern(f)}\de_x.$$
Bowen \cite{Bo1971, Bo1974} proved that the sequence above has a unique accumulation point $\mu$ which is $f$-invariant, ergodic, and is the unique measure maximizing the entropy, that is,
$$h_\mu(f)=h_{top}(f).$$

In this work, we give conditions on  $g\::\:X\to X$  an homeomorphism $C^0$-close to $f$  that ensure the existence and uniqueness of a probability measure maximizing the entropy. Under mild conditions such homeomorphism $g$ is semiconjugated to $f$, that is, there exists $h:X\to X$ continuous and onto such that $f\circ h= h\circ g,$  but  neither it is expansive  nor has the specification property.

Nevertheless, a weak form of expansivity is satisfied by $g$: if two points $x,y $ satisfies that
$$\dist(g^n(x),g^n(y))\leq \alpha, \, \quad\forall n\in\Z;$$
where $\al$ is the expansivity constant of $f$, then
$$\dist(g^n(x),g^n(y))\leq \eta, \,\quad \forall n\in\Z;$$
where $\eta$ is much smaller than $\al$ (and in fact depends on the $C^0$ distance between $f$ and $g$). In other words, $g$ losses the expansivity property and the specification property on a very tiny micro scale. Thus, if the perturbation does not increase the entropy, one may expect to have the same property as Bowen proved for $f$, since we can study the situation in a scale less than the expansivity constant of $f$ but larger compared to $\eta$. This is, morally, the heuristic reason for our result. See  Section~\ref{sec:statements} for detailed statements.

Our main motivation is to apply such conditions to the class of partially hyperbolic diffeomorphims, known as derived from Anosov. We conclude from our results that many such diffeomorphims admit a unique measure with the previous properties. See Section~\ref{sec:example} for more details.

We would like to point out that while we were working on this subject a preprint from J. Buzzi and T. Fisher (see \cite{BuFis2009}) came out. This work shows that a particular class of Ma\~{n}e's derived from Anosov diffeomorphims has a unique probability measure of maximal entropy  with methods similar to ours. Nevertheless, our method is a somewhat
more general and  it can be applied to more situations as we show in Section~\ref{sec:example} and to be extended to more general equilibrium states \cite{Samva2009}.

\indent{\bf Acknoledgements} We would like to thank Marcelo Viana for his comments and orientations. We also would like to thank Godofredo Iommi for enlighten us about the application to the suspention flows.

\section{Statements}\label{sec:statements}

Let $X$ be a compact metric space and $f\::X\to X$ be an $\al$-expansive homeomorphism with the specification property. For $n\geq 1$, denote by
$$\pern(f):=\{x\in X\::\: f^ n(x)=x\},$$
and define the sequence of $f$-invariants measures
\begin{equation}\label{eq:sucbowen}\mu_n:=\frac{1}{|\pern(f)|}\sum_{x\in \pern(f)}\de_x,\end{equation}
and define
\begin{equation}\label{eq:medbowen}
\mu:=\lim_{n\to \infty}\mu_n
\end{equation}

As we already said in the introduction, Bowen proved that the limit \eqref{eq:medbowen} in fact exists and $\mu$ is an $f$-invariant measure, ergodic,  and it is the unique measure maximizing the entropy of $f$.

Let $g\::\: X \to X$ be a homeomorphism semiconjugated to $f$, that means
there exists
$h:M\to M$ continuous and onto such that
$$f\circ h=h\circ g.$$

The semiconjugation $h$ allow to define the  equivalence relation:
$y\sim_h z$ if and only if $h(y)=h(z)$. For $x\in X$ we denote by
$[x]:=\{y\in X: h(y)=h(x)\}=h^{-1}(h(x))$ the equivalence class of $x$. We say that a class $[x]$ is {\em periodic} if $h(x)$ is a periodic point of $f$. A class $[x]$ is {\em trivial} if its cardinality is one. Otherwise, we say the class is non-trivial.

We recall that given a compact set $K$, the entropy of $g$ relative to $K$ is defined as follows:
for any $\epsilon >0$ and $n$, let $E_K(n,\epsilon)$ the maximum cardinality of a $(n,\epsilon)$-separated set (for $g$) and contained in $K$ and let
$$h_{top}(g,K)=\lim_{\epsilon\to 0}\limsup_{n\to\infty}\frac{1}{n}\log E_K(n,\epsilon).$$

Now, we are ready to state our conditions:

\begin{itemize}
\item[(H1)]  $h_{top}(g,[x])=0$ for any $x\in X.$
\item[(H2)] Each periodic class $[x]$ has a periodic point of $g$ with the same period of $h(x)$.
\item[(H3)] $\mu(\{h(x)\::\: [x] \mbox{ is trivial }\})=1.$

\end{itemize}

The condition (H2) allows to chose one periodic orbit for each periodic class, that is, we set
\begin{eqnarray*}
\perti(g):=&\{&x\in h^{-1}(\pern(f))\cap\pern(g)\::\: \\ &&\mbox{ if } o(x_1,g)\ne o(x_2,g) \mbox{ then } o(h(x_1),f)\ne o(h(x_2),f)\},
\end{eqnarray*}

and define
\begin{equation}\label{eq:sucsamva}\nu_n:=\frac{1}{|\perti(g)|}\sum_{x\in \perti(g)}\de_x,\end{equation}

\begin{maintheorem}\label{teo:main} If $f,g:X\to X$ are as above and assumptions (H1-H3) are satisfied then
\begin{equation}\label{eq:medsamva}
\nu=\lim_{n\to \infty}\nu_n
\end{equation} exists and it is $f$-invariant, ergodic,  and it is
the unique measure maximizing the entropy of $g$. In particular
$\nu$ equidistribute the periodic class of $g$.
\end{maintheorem}

\noindent{\bf Remark :} The condition (H2) can be removed. In fact, if $[x]$ is a periodic class with $h(x)\in\pern(f)$, then $[x]$ is compact and $g^m|[x]\::\:[x]\to[x]$ is an homeomorphism, where $m$ is the period of $h(x)$. So, there exists a $g^m$-invariant probability measure $\de_{[x]}$ suported on $[x]$. Of course $g_*\de_{[x]}$ is a $g^m$-invariant probability measure supported on $[h(x)]$ and $h_*\de_{[x]}=\de_{h(x)}$. In particular
$$\frac{1}{m}\sum_{k=0}^{m-1}g_*^k\delta_{[x]}$$
is a $g$-invariant probability measure supported on the orbit of the periodic class $[x].$

So we can define $\perti(g)$ as the set of equivalent classes that are fixed by $g^n$ and for each periodic class we select a $g$-invariant probability measure $\de_{[x]}$ as above (with $\de_{[g^k(x)]}=g^k_*\de_{[x]})$.  Now \eqref{eq:sucsamva} can be written as
\begin{equation}\label{eq:sucsamva2}\nu_n:=\frac{1}{|\perti(g)|}\sum_{{[x]}\in \perti(g)}\de_{[x]}.\end{equation}
Note that (H2) imples that $\de_{[x]}$ can be chosed a Dirac measure supported on a periodic point of $[x]$. Nevertheless, as we will see in the proof of Lemma~\ref{le:uno}, it is enough that there exist a biyection between $\pern(f)$ and $\perti(g)$ and
$$h_*\de_{[x]}=\de_{h(x)}.$$

 \section{Proof of Theorem~\ref{teo:main}}\label{sec:proof}

Let $\nu$ be any accumulation point of the sequence
$\nu_n.$ We will prove that $\nu$ is the unique measure of
maximal entropy, and hence $\nu$ will be the limit of $\nu_n$ and the result will follows.

\begin{lema}\label{le:uno}
It holds that $h_*\nu=\mu.$
\end{lema}
\begin{proof}
First, note that $h_*\nu_n=\mu_n$, for all $n\geq 1$. In fact, for
every $A\subseteq X$ Borelean,
 if $x\in X$ is any point  then
$$\de_x(h^{-1}(A))=\de_{h(x)}(A).$$
Taking this last equality and (H2) into consideration, we conclude that
\begin{eqnarray*}
h_*\nu_n(A)&=&\nu_n(h^{-1}(A))\\
&=&\frac{1}{|\perti(g)|}\sum_{x\in \perti(g)}\de_x(h^{-1}(A))\\
&=&\frac{1}{|\pern(f)|}\sum_{h(x)\in \pern(f)}\de_{h(x)}(A)\\
&=&\mu_n(A).
\end{eqnarray*}
From the continuity of $h_*$ we have that
$$h_*\nu=\mu.$$
\end{proof}
\begin{lema} The measure $\nu$ is of maximal entropy, that is,
$$h_\nu(g)=h_{top}(g).$$
\end{lema}
\proof Let $\mathcal{P}$ be a partition of $X$ and
$\tilde{\mathcal{P}}=h^{-1}(\mathcal{P})$. Then
$h_\nu(g,\tilde{\mathcal{P}})=h_\mu(f,\mathcal{P})$ and so
$h_{\mu}(f)\leq h_\nu(g)$. On the other hand,  Bowen formula \cite{Bo1971b} states
that
$$h_{top}(g)\leq h_{top}(f)+ \sup_{x\in X}h_{top}(g,[x]).$$
Therefore, from (H1) and the variational principle we conclude that
$$h_{top}(g)\leq h_{top}(f)=h_\mu(f)\leq h_\nu(g)\leq h_{top}(g).$$
\endproof
We say that $A$ is saturated if $A=h^{-1}(h(A))$. In general,  the {\em saturation} of $A\subseteq X$ is defined as ${\rm sat}(A):=h^{-1}(h(A))$.
 Note that $\nu({\rm sat}(A))=\mu(h(A))$.
\begin{lema}\label{le:saturado} For ever Borel set $A$ we have
 $\nu(A)=\nu({\rm sat}(A))$.
\end{lema}
\proof Let $\tilde{X}=\{x\in X\::\: [x]=\{x\}\}$. From (H3) and the
fact that $h_*\nu=\mu$ we have that $\nu(\tilde{X})=1$. For
$A\subseteq X$ Borelean, we have
$$\nu({\rm sat}(A))=\nu({\rm sat}(A)\cap\tilde{X})=\nu(A\cap\tilde{X})=\nu(A).$$
\endproof

\begin{cor}
The probability measure $\nu$ is ergodic.
\end{cor}
\begin{proof}
From Lemma~\ref{le:saturado} follows that if $P$ is a $g$-invariant
subset, then
$$\nu(P)=\nu(h^{-1}(h(P))=\mu(h(P)).$$
Since $h(P)$ is $f$-invariant and $\mu$ is ergodic, then $\nu$ is
ergodic.
\end{proof} .

\begin{lema}
Let $\eta$ be a $g$-invariant probability measure and assume that $\eta$ is singular with respect to $\nu.$ Then
$$h_\eta(g)<h_{top}(g)$$
\end{lema}

\begin{proof}
Let $\rho=h_*\eta.$ It follows that $\rho$ is singular with respect to $\mu.$  The Ledrappier-Walter's formula \cite{LeWal1977} states that
$$h_\eta(g)\leq h_\rho(f) + \int_Xh_{top}(g,h^{-1}(x))d\rho(x).$$
\noindent and from (H1) it follows that
$$h_{\eta}(g)\leq h_{\rho}(f).$$
Bowen proved \cite{Bo1974} that $h_{\rho}(f)<h_{top}(f)=h_{top}(g)$ and the result follows.
\end{proof}

Now, we are about to finish the proof of Theorem~\ref{teo:main}. Let $\eta$ be any $g$-invariant probability measure such that $h_{\eta}(g)=h_{top}(g).$ We can write $\eta=\alpha\eta_1+(1-\alpha)\eta_2$ for some $\alpha\in [0,1]$ such that $\eta_i$ are probablity measures, $\eta_1<<\nu$ and $\eta_2$ is singular with respect to $\nu.$ It follows that
$$h_{top}(g)=h_\eta(g)=\alpha h_{\eta_1}(g)+(1-\alpha)h_{\eta_2}(g)\leq h_{top}(g).$$
The previous lemma implies that $\alpha=1$, that is, $\eta$ is absoluty continuous with respect to $\nu.$ As $\nu$ is ergodic we have that $\eta=\nu.$ This completes the proof of the theorem.

\section{Examples}\label{sec:example}

This section is devoted to show some examples where our theorem applies.

\subsubsection{Example 1: Blowing up periodic point.}

The first and almost trivial example is when we blow up a fixed point. Let $f:M\to M$ be an expansive homeomorphism of the manifold $M$ with the specification property and let $p$ be a fixed point. Let $B(p)$ be a closed ball around $p$ and let $h:M\to M$ be a continuos map such that $h(B)=p$ and $h|(M\setminus B):M\setminus B\to M\setminus\{p\}$ is a bijection. Let $g_B:B\to B$ be a homeomorphism such that it has zero entropy on $B$ and assume that   $g:M\to M$ defined as $g(x)=g_B(x)$ if $x\in B$ and $g(x)=h^{-1}(f(h(x))$ otherwise is a homemorphism. Thus, $g$ satisfies our theorem and hence has a unique measure of maximal entropy.

\subsection{Derived from Anosov examples}

For the sake of completeness we include  known results about $C^0$ perturbation of Linear Anosov diffeomorphisms.
Let $A:\R^n\to\R^n$ be a hyperbolic linear map $\R^n=E^s\oplus E^u.$ We take an adapted norm on each subspace:
$$\|A|E^s\|_s<a<1,\,\,\,\,\|A^{-1}|E^u\|_u<a<1$$
and in $\R^n$ we take the maximum norm with respect to the decomposition $E^s\oplus E^u$, that is $\|x\|=\|(x^s,x^u)\|=\max\{\|x^s\|_s,\|x^u\|_u\}.$ Let $r>0$. We say that $\{x_n:n\in \Z\}$ is a $r$-pseudo orbit (for $A$) if $\|Ax_n-x_{n+1}\|\le r$ for all $n\in \Z.$

\begin{lema}
Let $\{x_n\}_{n\in\Z}$ be an $r$-pseudo orbit. Then it is $\frac{r}{1-a}$ shadowed by a true orbit, i. e., there exists $y$ such that, for every $n\in \Z$, then
$$\|A^ny-x_n\|\le \frac{r}{1-a}.$$
\end{lema}
\begin{proof}
Let $x_n=(x_n^s,x_n^u).$ It is not dificult to see (by the uniform contraction in $E^s$) that, for every $n\geq 0$ and all $m\in\Z$, then
$$\|A^nx^s_m-x^s_{n+m}\|\leq \frac{r}{1-a}.$$
Thus, if $y^s$ is an acumulation point of $A^nx^s_{-n}$ then, for every $n\in\Z$, $\|A^ny^s-x^s_n\|\le \frac{r}{1-a}.$
An anologous argument show that there is a point $y^u$ such that for every $n\in\Z$, $\|A^ny^u-x^u_n\|\le \frac{r}{1-a}.$ Therefore, the point $y=(y^s,y^u)$ is the desired point. Uniqueness follows since $A$ is infinitely expansive.
\end{proof}

\begin{cor}
Let $A$ be hyperbolic and let $G:\R^n\to\R^n$ be a homeomorphism such that $\|Ax-G(x)\|\leq r$ for all $x\in\R^n.$ Then there exists $H:\R^n\to\R^n$ continuous and onto such that $A\circ H=H\circ G.$ Moreover $\|H(x)-x\|<\frac{r}{1-a}.$
\end{cor}
\begin{proof}
Any $G$-orbit is an $r$ pseudo orbit for $A.$ Therefore, just define $H(x)$ as the unique point such that
$$\|A^n(H(x))-G^n(x)\|\leq\frac{r}{1-a}$$
and the result follows.
\end{proof}

\begin{cor}
Let $f:\T^n\to\T^n$ be a linear Anosov map. Then there exists $C>0$ such that for any small $r$ and any $g:\T^n\to\T^n$ with $\dist_{C^0}(f,g)<r$ there exists $h:\T^n\to\T^n$ continuous and onto, $\dist_{C^0}(h,Id)<Cr,$ and $f\circ h= h\circ g.$ Furthermore, if $\alpha$ is a expansivity constant of $f$ and $Cr<\alpha/4$ , then if $\dist(g^n(x),g^n(y))\leq \alpha/2$ for all $n\in \Z$, then $h(x)=h(y).$
\end{cor}

\begin{proof}
Let $A\in SL(n,\Z)$ be the lift of $f$ to $\R^n$ and let $G$ be a lift of $g$ such that $\|Ax-G(x)\|<r$ for all $x.$ If $r$ is small, then $G$ is isotopic to $A$ and hence $G=A+p$ where $p(x+\Z^n)=p(x).$ Notice that $G^n=A^n+p_n$ where $p_n(x+\Z^n)=p_n(x).$
Let $H$ from the previous corollary and let $t\in\Z^n.$ Then:
\begin{eqnarray*}
\|A^n(H(x)+t)-G^n(x+t)\|&=&\|A^n(H(x))+A^nt-G^n(x)-A^nt\|\\ &=&\|A^n(H(x))-G^n(x)\|\leq\frac{r}{1-a},
\end{eqnarray*}
for all $n\in\Z$ and hence $H(x+t)=H(x)+t.$ Thus $H$ defines the desired $h$ on $\T^n.$
The  last part follows since $\dist(f^n(h(x)),f^n(h(y)))<\alpha$ for all $n\in\Z.$
\end{proof}

\subsubsection{Example 2: Ma\~{n}e's derived from Anosov}

Let $\T^n$, $n\geq 3$, be the torus $n$-dimentional and $f\::\:\T^ n\to \T^n$ be a (linear) Anosov diffeomorphism with expansivity constant $\alpha$. Assume that the tangent bundle of $\T^n$ admits the $Df$-invariant splitting
$$T\T^n=E^{ss}\oplus E^u \oplus E^{uu},$$
with $\dim E^u=1$ and $$\la_s:=\|Df|E^{ss}\|\,\,\,\,\,\la_u:=\|Df|E^u\|\,\,\,\,\la_{uu}:=\|Df^{-1}|E^{uu}\|$$ satisfying the relation $\la_s<1<\la_u<\la_{uu}$. Let $p$ be a fixed point of $f.$ Let $r>0$ small enough (to be fixed later). Consider the open ball $B(p,r)$. Then deform the Anosov diffeomorphim $f$ inside $B(p,r)$ passing through a flip bifurcation along the central unstable foliation $\fu(p)$ and then we obtain three fixed point, two of them with stability index equal to $\dim E^s$ and the other one with stability index equal to $\dim E^s+1$ as it is  shown in Figure~\ref{fig1}.

\begin{figure}[h]
\begin{flushleft}
\vspace*{2.5cm}
\hspace*{3cm}
\psline[linewidth=0.5pt]{->}(-1,0)
\psline[linewidth=0.5pt]{->}(1,0)
\psline[linewidth=0.5pt]{-}(-1,0)(-2,0)
\psline[linewidth=0.5pt]{-}(1,0)(2,0)
\psline[linewidth=0.5pt]{->}(0,-1)
\psline[linewidth=0.5pt]{->}(0,1)
\psline[linewidth=0.5pt]{>-}(0,-1)(0,-2)
\psline[linewidth=0.5pt]{>-}(0,1)(0,2)
\psline[linewidth=0.5pt]{-<}(-0.7,-0.7)
\psline[linewidth=0.5pt]{->}(-1.5,-1.5)(-0.7,-0.7)
\pscircle[linestyle=dashed,dash=3pt 2pt,linewidth=0.5pt]{1.3}
\rput(.2,.4){$p$}
\rput(1.3,1.3){$f$}
\rput(1.3,-1.4){$B(p,r)$}
\psdot*[dotscale=.7](0,0)
\psline[linewidth=0.5pt]{->}(6,0)(5,0)
\psline[linewidth=0.5pt]{->}(6,0)(7,0)
\psline[linewidth=0.5pt]{-}(5,0)(4,0)
\psline[linewidth=0.5pt]{-}(7,0)(8,0)
\psline[linewidth=0.5pt]{->}(6,0)(6,-1)
\psline[linewidth=0.5pt]{->}(6,0)(6,1)
\psline[linewidth=0.5pt]{>-}(6,-1)(6,-2)
\psline[linewidth=0.5pt]{>-}(6,1)(6,2)
\psline[linewidth=0.5pt]{-<}(6,0)(5.3,-0.7)
\psline[linewidth=0.5pt]{->}(4.5,-1.5)(5.3,-0.7)
\pscircle[linestyle=dashed,dash=3pt 2pt,linewidth=0.5pt](6,0){1.3}
\rput(6.2,.4){$p$}
\rput(7.3,1.3){$g$}
\rput(7.3,-1.4){$B(p,r)$}
\psline[linewidth=0.5pt]{->}(5.5,0)(5.5,-1)
\psline[linewidth=0.5pt]{->}(5.5,0)(5.5,1)
\psline[linewidth=0.5pt]{>-}(5.5,-1)(5.5,-2)
\psline[linewidth=0.5pt]{>-}(5.5,1)(5.5,2)
\psline[linewidth=0.5pt]{->}(6.5,0)(6.5,-1)
\psline[linewidth=0.5pt]{->}(6.5,0)(6.5,1)
\psline[linewidth=0.5pt]{>-}(6.5,-1)(6.5,-2)
\psline[linewidth=0.5pt]{>-}(6.5,1)(6.5,2)
\psline[linewidth=0.5pt]{-<}(6.5,0)(5.8,-0.7)
\psline[linewidth=0.5pt]{->}(5,-1.5)(5.8,-0.7)
\psline[linewidth=0.5pt]{-<}(5.5,0)(4.8,-0.7)
\psline[linewidth=0.5pt]{->}(4,-1.5)(4.8,-0.7)
\psline[linewidth=0.5pt]{->}(6.5,0)(6.2,0)
\psline[linewidth=0.5pt]{->}(5.5,0)(5.8,0)
\psdot*[dotscale=.7](6,0)
\psdot*[dotscale=.7](6.5,0)
\psdot*[dotscale=.7](5.5,0)
\vspace*{2.5cm}
\caption{}\label{fig1}
\end{flushleft}
\end{figure}
We call $g$ the diffeomorphism obtained in this way. It made be done in such a way  that $g$ is partially hyperbolic with subbundless $E^{ss}\oplus E^c\oplus E^{uu}$ with $\dim E^c=1$ and such that  $\dist_{C^0}(f,g)<r$. Then, $g$ is semiconjugated to $f$ by $h:M\to M$ which is close to the identinty of the same order as $r,$ that is $\dist(h(x),x)\le Cr$ for some constant $C.$

Notice that, since $E^c$ is one dimensional (and continuous) it is integrable (not necessarily unique). And hence, $E^{ss}\oplus E^c$ and $E^c\oplus E^{uu}$ are integrable too. We write $W^{cs}_\epsilon$ the ball of radius $\epsilon$ in $W^{cs}(x)$ (for some central stable plaque through $x$). Once we fixed a central stable manifold through $x$, say $W^{cs}(x)$, we set $W^{cs}(g^n(x)):=g^n (W^{cs}(x)).$

We will assume the following properties of $g:$ (see Figure \ref{fig2} and \ref{fig3})

\begin{itemize}
\item We may choose, $\delta$ and $\epsilon<\alpha/2$ such that if $d(x,y)<\delta$ then $W^{cs}_\epsilon(x)\cap W^{uu}_\epsilon(y)$ is nonempty and consists of a single point, say $z.$ Moreover, $z\in W^{uu}_{\epsilon/\la_{uu}}(g(y))$ and $g(z)\in W^{cs}_\epsilon(g(x)).$
\end{itemize}
We remark that $\delta$ and $\epsilon$ are independent of $r$ small, they depend on the $C^0$ closeness of $E^{ss}(f), E^c(f), E^{uu}(f)$ with the respective bundles for $g$.

\begin{figure}[h]
\begin{flushleft}
\vspace*{4cm}
\hspace*{-1.5cm}
\rput(5.7,2.5){$W_\epsilon^{cs}(x)$}
\rput(3.7,3.5){$W_\epsilon^{uu}(y)$}
\pscurve[linewidth=0.5pt,showpoints=false]{-}(3,2)(4.5,2.2)(5,2.3)
\pscurve[linewidth=0.5pt,showpoints=false]{-}(1,0)(2.5,1.7)(3,2)
\pscurve[linewidth=0.5pt,showpoints=false]{-}(1,0)(2.5,.2)(4,.3)
\pscurve[linewidth=0.5pt,showpoints=false]{-}(4,0.3)(4.5,.8)(5,2.3)
\psdot*[dotscale=.5](3.3,1.1)
\rput(3.5,1.05){$z$}
\psdot*[dotscale=.5](2.3,.8)
\rput(2.5,.8){$x$}
\psdot*[dotscale=.5](3.3,2.8)
\rput(3.5,2.8){$y$}
\psline[linewidth=0.5pt]{->>}(3.3,1.2)(3.3,1.8)
\psline[linewidth=0.5pt]{-}(3.3,1.8)(3.3,3.3)
\psline[linestyle=dashed,dash=3pt 2pt,linewidth=0.5pt]{-}(3.3,1.2)(3.3,.3)
\psline[linewidth=0.5pt]{->>}(3.3,.2)(3.3,-.3)
\psline[linewidth=0.5pt]{-}(3.3,-.3)(3.3,-.8)
\psline[linewidth=0.5pt]{->}(6.5,2)(8.5,2)
\rput(7.5,2.3){$g$}
\rput(12.7,2.5){$W_\epsilon^{cs}(g(x))$}
\rput(10.7,3.5){$W_{\epsilon/\lambda_{uu}}^{uu}(g(y))$}
\pscurve[linewidth=0.5pt,showpoints=false]{-}(10,2)(11.5,2.2)(12,2.3)
\pscurve[linewidth=0.5pt,showpoints=false]{-}(8,0)(8.3,.2)(9.5,1.7)(10,2)
\pscurve[linewidth=0.5pt,showpoints=false]{-}(8,0)(9.5,.2)(11,.3)
\pscurve[linewidth=0.5pt,showpoints=false]{-}(11,.3)(11.5,.8)(11.7,2)(12,2.3)
\psdot*[dotscale=.5](10.3,1.1)
\rput(11,1.05){$g(z)$}
\psdot*[dotscale=.5](9.3,.8)
\rput(9.8,.8){$g(x)$}
\psdot*[dotscale=.5](10.3,2.8)
\rput(10.8,2.8){$g(y)$}
\psline[linewidth=0.5pt]{->>}(10.3,1.2)(10.3,1.8)
\psline[linewidth=0.5pt]{-}(10.3,1.8)(10.3,3.3)
\psline[linestyle=dashed,dash=3pt 2pt,linewidth=0.5pt]{-}(10.3,1.2)(10.3,.3)
\psline[linewidth=0.5pt]{->>}(10.3,.2)(10.3,-.3)
\psline[linewidth=0.5pt]{-}(10.3,-.3)(10.3,-.8)
\caption{}\label{fig2}
\end{flushleft}
\end{figure}

\begin{itemize}
\item We also assume that $[x,y]^c$ is a central segment whose length is less than $\epsilon$ and $\dist(g(x),g(y))<\delta$ then $g([x,y]^c)$ is a central segment with length less than $\epsilon.$ The same with $g^{-1}.$
\end{itemize}
\begin{figure}[h]
\begin{flushleft}
\vspace*{2cm}
\hspace*{.5cm}
\psline[linewidth=0.4pt,linearc=0.6]{->}(5,1.3)(5.5,1.5)(6,1.3)
\rput(5.5,1.8){$g$}
\pscurve[linewidth=0.5pt,showpoints=false]{-}(1,0)(2,0.3)(3,1)(3.5,1.4)(4,2)
\pscurve[linewidth=1pt,showpoints=false]{-}(2,0.3)(3,1)(3.5,1.4)
\psdot*[dotscale=.7](2,0.3)
\rput(2.3,0.2){$x$}
\psdot*[dotscale=.7](3.5,1.4)
\rput(3.8,1.3){$y$}
\rput(2.5,1.3){$[x,y]^c$}
\rput(1,-.5){$\mathcal{F}^c(x)$}
\pscurve[linewidth=0.5pt,showpoints=false]{-}(6.5,0)(7.5,0.6)(8.5,.9)(9,1.3)(10,1.5)
\pscurve[linewidth=1pt,showpoints=false]{-}(7.5,0.6)(8.5,.9)(9,1.3)
\rput(7.8,1.3){$g([x,y]^c)$}
\rput(6.5,-.5){$\mathcal{F}^c(g(x))$}
\psdot*[dotscale=.7](7.5,.6)
\rput(7.8,0.35){$g(x)$}
\psdot*[dotscale=.7](9,1.3)
\rput(9.3,1.05){$g(y)$}
\vspace*{.5cm}
\caption{}\label{fig3}
\end{flushleft}
\end{figure}

\begin{lema}
Assume that $r$ is such that $Cr<\min\{\delta,\alpha/4\}$. Let $g:\T^n\to\T^n$ be such that $\dist_{C^0}(f,g)<r$ and $g$ satisfies the above assumptions. Let $h$ be the semiconjugacy $Cr$ close to the identity between $f$ and $g.$ Then, an equivalent class consists of a single (closed) central segment (possibly trivial) of length less than $2Cr.$
\end{lema}

\begin{proof}
Let $x,y$ be such that $h(x)=h(y).$ Fix a central stable plaque through $x$ and let $z=W_\epsilon^{cs}(x)\cap W_\epsilon^{uu}(y).$ Then, since $\dist(g(x),g(y))<\delta$ it follows that
$$W^{cs}_\epsilon(g(x))\cap W^{uu}_\epsilon(g(y))=\{g(z)\}.$$
By an induction argument we conclude that $g^n(z)\in W^{uu}_\epsilon(g^n(y)).$ Thus, $z=y$ and so $y\in W^{cs}_\epsilon(x).$
With a very similar argument, arguing with $W^c_\epsilon(x)$ and $W^{ss}(y)$ inside $W^{cs}_\epsilon(x)$ we conclude that $y\in W^c(x).$  Let $[x,y]^c$ be the central segment in the chosen $W^c(x)$ between $x$ and $y$, it has length less than $\epsilon.$ Now we conclude that $g^n([x,y]^c)$ has length less than $\epsilon$ for all $n.$ This implies that $h([x,y]^c)=h(x)=h(y).$

\end{proof}

\begin{cor}
Let $g$ be as above. Then {\rm (H1)} and {\rm (H2)} of Theorem \ref{teo:main} hold.
\end{cor}
\begin{proof}
Condition (H2) holds since a periodic segment has a periodic point of the same period. Condition (H1) holds since, for every $\epsilon$ and $n$ the cardinality of a maximal $(n,\epsilon)$ separated set in an equivalent class is bounded by $n2Cr/\epsilon.$
\end{proof}

In order to verify that hypothesi (H3) is satisfied, consider the following:

\begin{enumerate}
\item Set $\sigma=\inf\{\|Dg(x)|E^c(x)\|\::\:x\in M\}.$

\item Let $m=\min\{j\ge 1\::\: \sigma \la_u^j>1\}.$

\item Let $\rho$ be such that $\mu(M\setminus B(p,\rho))\ge 1-\frac{1}{2m}$ where $\mu$ is the maximal entropy
probability measure of $f.$

\item Let $r$ be such that $2Cr<\rho/2.$
\end{enumerate}

We will show that these conditions imply (H3) holds. We remark under the above setting that $g$ could be very far in the $C^1$ topology from any Anosov diffeomorphism (if $\sigma$ is very small) but in that case the $C^0$ distance between $f$ and $g$ should be very small.

\begin{lema}
Under the above setting, condition {\rm (H3)} is satisfied, that is, if $x$ is $\mu$-generic, then $h^{-1}(x)$ is a single point.
\end{lema}

\begin{proof}
Let $x$ be a $\mu$- generic point.
  Then Birkhoff's theorem ensure that
$$\lim_{n\to\infty}\frac1n\{j\::\:f^j\notin B(p,\rho)\}\geq 1-\frac{1}{2m},$$
and so, for $n\geq n_0$, with $n_0$ large we have
$$\frac1n\{j\::\:f^j(x)\notin B(p,\rho)\}\geq1-\frac1m.$$
Now, take a point $y\in h^{-1}(x)$. Then, for every $n\in\Z$ we have
$$\dist(g^ n(y),f^n(x))<Cr.$$
So, for $n\geq n_0$,
$$\frac1n\{j\::\:g^j(y)\notin B(p,\rho/2)\}\geq 1-\frac1m.$$
We claim that $h^{-1}(x)$ consists of one point.  Otherwise, let $y_1, y_2\in h^ {-1}(x)$. Take $n\geq n_0$ such that
$$(\si\la_u^m)^{-n}2Cr<\dist(y_1,y_2).$$
Since
$$\dist(g^n(y_1), g^n(y_2))<2Cr$$ and the central segment between $y_1,y_2$ is also in the equivalent class
 we have that
\begin{eqnarray*}
\dist(y_1,y_2)& =& \dist(g^{-n}(g^n(y_1)),g^{-n}(g^n(y_2)))\\ &\leq& (\si\la_u^m)^{-n}\dist(g^n(y_1), g^n(y_2))\\ &<&(\si\la_u^m)^{-n}2Cr<\dist(y_1,y_2),
\end{eqnarray*}

\noindent a contradiction. Therefore, for any $\mu$-generic point $x$, we have that $h^{-1}(x)$ consist of a single point and (H3) is satisfied.

\end{proof}

\begin{cor}
Let $g$ as above. Then $g$ has a unique probability measure of maximal entropy.
\end{cor}

\subsubsection{Example 3: Mixed Ma\~{n}e's derived from Anosov}

Let $f:\T^n\to \T^n$ be a (linear) Anosov diffeomorphism, $n\ge 4$ such that $T\T^n=E^{ss}\oplus E^s\oplus E^u\oplus E^{uu}$ with $\dim E^s=\dim E^u=1$ with rate of contraction/expansion as $\la_{ss}<\la_s<1<\la_u<\la_{uu}.$ Taking a power of $f$ if necessary, assume that $f$ has two different fixed points $p$ and $q$. Let $r>0$ be small (to be fixed later) and deform $f$ inside $B(p,r)$ and $B(q,r)$ like the Ma\~{n}e's derived from Anosov: in $B(p,r)$ we perform a flip perturbation along $E^s$ and on $B(q,r)$ we perform a flip bifurcation along $E^u.$ In this way we can obtain $g$ satisfying the following:
\begin{itemize}
\item $g$ is partially hyperbolic: $T\T^n=E^{ss}\oplus E^{cs}\oplus E^{cu}\oplus E^{uu}$ which is dominated (each subbundle dominates the previous ones by a factor $a<1$), and $\dim E^{cs}=\dim E^{cu}=1.$ These subbundles are $C^0$ close to the respective ones of $f.$
\item $\dist_{C^0}(f,g)<r.$
\item If $\dist(x,y)<2Cr$ then $\frac{\|Dg|E^{ci}(x)\|}{\|Dg|E^{ci}(y)\|}<a^{-1/4},\,i=s,u$
\item $Df|E^{cs}(x)$ is uniformly contracting outside $B(p,r)$ with rate $\la_s.$
\item $Df|E^{cu}(x)$ is uniformly expanding  outside $B(q,r)$ with rate $\la_u.$
\item Let $\si_1=\sup\{\|Dg|E^{cs}(x)\|\::\:x\in M\}$ and $\si_2=\inf\{\|Dg|E^{cu}(x)\|\::\:x\in M\}$ and $m$ be such that $\si_1\la_s^m<1$ and $\si_2\la_u^m>1.$ Let $\rho$ such that $\mu(M\setminus B(j,\rho))\ge 1-1/2m, \,j=p,q$ where $\mu$ is the Bowen measure of $f.$ We assume that $2Cr<\rho/2.$

\end{itemize}

We show that $g$ falls into the assumptions of our Theorem and hence has a unique measure of maximizing entropy. This example share similarities with the previous example and so we will give an outline of the proof. Let $h:\T^n\to\T^n$ the semicongugacy. We will assume for the sake of simplicity that $E^c=E^{cs}\oplus E^{cu}$ is uniquely integrable. Notice that $E^{cs}$ and $E^{cu}$ are one-dimensional and hence they are integrable. Let $J$ be a segment tangent to $E^{cs}$, we say that $E^{cs}$ is uniquely integrable through $J$ if any maximal integral curve of $E^{cs}$ through any point of $J$ must contain $J.$ Analogously for $E^{cu}.$

\begin{lema}
Let $x\in \T^n$ be any point. Then, one and only one of the following hold:
\begin{enumerate}
\item $h^{-1}(x)$ consists of a single point.
\item $h^{-1}(x)$ is a segment tangent to $E^{cs}$ of length less than $2Cr.$
\item $h^{-1}(x)$ is a segment tangent to $E^{cu}$ of length less than $2Cr.$
\item $h^{-1}(x)$ is a square tangent to $E^{cs}\oplus E^{cu}$:
 \begin{itemize}
\item for each $y\in h^{-1}(x)$ we have that $W^{cs}_\epsilon(y)\cap h^{-1}(x)$ is a central stable segment that we denote by $J^{cs}(y)$ and $E^{cs}$ is uniquely integrable through $J^{cs}(y).$ Similar for $E^{cu}.$
 \item If $y$ and $z$ are in $h^{-1}(x)$ then, $\emptyset\neq J^{cs}(y)\cap J^{cu}(z) \in h^{-1}(x).$
 \end{itemize}
\end{enumerate}
\end{lema}
\begin{proof}
Assume that $h^{-1}(x)$ is not trivial, and let $y,z\in h^{-1}(x)$ be two different points.
By a similar argument as the previous example we conclude that $y\in W^c_\epsilon(z).$ And also, if $z\in W^{cs}(y)$ then $[y,z]^{cs}\subset h^{-1}(x).$ This means that $W^{cs}_\epsilon (y)\cap h^{-1}(x)$  is a segment, say $J^{cs}(y)$ whose length remains bounded in the future and in the past and, by the domination in $E^{cs}\oplus E^{cu}$ we conclude that  $E^{cs}$ is uniquely integrable through $J^{cs}(y).$  Similar if $z \in W_\epsilon^{cu}(y).$

Assume also that neither $2)$ nor $3)$ hold. Consider local central integral curves $W_\epsilon^{cs}(y)$ and $W^{cu}_\epsilon(z)$ and call $w$ the point of intersection. Although they may not have rate of expansion or contraction, a similar argument can be done so that $h(w)=h(z)=h(y).$ Therefore, $\{w\}= J^{cs}(y)\cap J^{cu}(z) \in h^{-1}(x).$
\end{proof}
\begin{figure}[ht]
\begin{flushleft}
\vspace*{4cm}
\hspace*{2.5cm}
\pscurve[linewidth=0.5pt,showpoints=false]{-}(0,0)(.3,1)(.5,1.5)(0.7,3)
\pscurve[linewidth=0.5pt,showpoints=false]{-}(5,0)(5.3,1)(5.5,1.5)(5.7,3)
\pscurve[linewidth=0.5pt,showpoints=false]{-}(0.7,3)(1,3.1)(3.5,3.5)(5.7,3)
\pscurve[linewidth=0.5pt,showpoints=false]{-}(0,0)(1,.1)(3.5,.5)(5,0)
\rput(6.2,3.2){$\mathcal{F}^c(x)$}
\pscurve[linewidth=0.5pt,showpoints=false]{-}(1,.5)(1.2,1.3)(1.5,2.1)(1.6,2.7)
\pscurve[linewidth=0.5pt,showpoints=false]{-}(2,.8)(2.2,1.6)(2.5,2.4)(2.6,3.1)
\pscurve[linewidth=0.5pt,showpoints=false]{-}(3,1)(3.2,1.8)(3.5,2.6)(3.6,3.3)
\pscurve[linewidth=0.5pt,showpoints=false]{-}(4,.8)(4.2,1.6)(4.5,2.4)(4.6,3.1)
\pscurve[linewidth=0.5pt,showpoints=false]{-}(0.8,0.7)(1.8,1.2)(3.5,1.4)(4.8,1)
\pscurve[linewidth=0.5pt,showpoints=false]{-}(0.8,1.2)(1.8,1.7)(3.5,1.9)(4.8,1.7)
\pscurve[linewidth=0.5pt,showpoints=false]{-}(0.8,1.6)(1.8,2.1)(3.5,2.3)(4.9,2.1)
\pscurve[linewidth=0.5pt,showpoints=false]{-}(1,2)(1.8,2.5)(3.5,2.8)(4.9,2.5)
\psdot*[dotscale=.7](1.09,.85)
\rput(1.3,.5){$y$}
\psdot*[dotscale=.7](3.4,2.3)
\rput(3.6,2.1){$z$}
\psdot*[dotscale=.7](1.45,1.97)
\rput(1.8,2.3){$w$}
\pscurve[linewidth=0.5pt,showpoints=false]{->}(1.6,3.9)(1.5,3.4)(1.6,2.9)
\rput(1.5,4.1){$W_\epsilon^{cs}(y)$}
\pscurve[linewidth=0.5pt,showpoints=false]{->}(6.5,2.2)(6,2.25)(5,2.2)
\rput(7.3,2.2){$W^{cu}_\epsilon(z)$}
\caption{}\label{fig4}
\end{flushleft}
\end{figure}

\begin{cor}
Conditions {\rm (H1)} and {\rm (H2)} are satisfied.
\end{cor}
\begin{proof}
We need only to check (H1) in case (iv) above. We the above notations, we observe that $g(J^{cs}(y))=J^{cs}(g(y))$ and similar for $cu.$ Therefore, the product structure is invariant and it is not difficult to see that the maximal cardinality of a $(n,\epsilon)$ separated set in the equivalent class has at most polynomial growth. Anyway, we can also say a bit more on the structure of the iterates of an equivalent class.  Let $y$ be in the class and consider $J^{cu}(y).$ Since the length remains bounded in the future, we claim that $\|Dg^n|E^cu(y)\|\le a^{-1/2}$ for any $n$ large enough. Otherwise, since for any $w\in J^{cu}(y)$ we have that
 $$\|Dg^n|E^{cu}(w)\|\ge \|Dg^n|E^{cu}(y)\|a^{1/4}\ge a^{-1/4}$$
and therefore the segment $J^{cu}(y)$ will be, for $n$ large, larger than $2Cr,$ a contradiction, and the claim is proved. Now,
by domination, we conclude that $\|Dg^n|E^{cs}(y)\|\le a^{1/2}$ for every $n$ large enough, and in fact the same holds for any $w\in J^{cu}(y).$ Hence,  the length of $J^{cs}(w)$ will decrease exponentially fast. We remark that a similar property holds in the past: the $cu$ segments are contracted exponentially fast.  Then, it is not difficult to see that for any $\epsilon>0$ the cardinality of a maximal $(n,\epsilon)$ separated set in the class of $y$ is the same as the cardinality of a maximal $(n,\epsilon)$ separated set in $J^{cu}(y)$ for $n$ large enough. Thus, $h_{top}(g,[y])=0.$
Condition (H2) holds trivially in the first three cases, and also in the fourth one by use, for instance, of Brouwer's fixed point theorem. Nevertheless, by the structure given above, case (iv) can not be a periodic class.
\end{proof}

\begin{figure}[ht]
\begin{flushleft}
\vspace*{4cm}
\hspace*{-1cm}
\psline[linewidth=0.5pt]{->}(0.7,.5)(1,2.2)
\rput(.5,1.4){$cs$}
\pscurve[linewidth=0.5pt,showpoints=false]{-}(1,.5)(1.3,1)(1.7,3.3)(2,4)
\pscurve[linewidth=0.5pt,showpoints=false]{-}(1.2,.5)(1.5,1)(1.9,3.3)(2.2,4)
\pscurve[linewidth=0.5pt,showpoints=false]{-}(1.4,.5)(1.7,1)(2.1,3.3)(2.4,4)
\pscurve[linewidth=0.5pt,showpoints=false]{-}(1.6,.5)(1.9,1)(2.3,3.3)(2.6,4)
\psline[linewidth=0.5pt]{-}(1,.5)(1.6,.5)
\psline[linewidth=0.5pt]{-}(1.45,1.5)(2.05,1.5)
\psline[linewidth=0.5pt]{-}(1.6,2.7)(2.2,2.7)
\psline[linewidth=0.5pt]{-}(2,4)(2.6,4)
\psline[linewidth=0.5pt]{->}(1,0)(2.5,0)
\rput(1.7,-.5){$cu$}
\psline[linewidth=0.5pt]{->}(4.7,.5)(5,2.2)
\rput(4.5,1.4){$cs$}
\psline[linewidth=0.5pt]{->}(5,0)(6.5,0)
\rput(5.7,-.5){$cu$}
\rput(7,-1){$h^{-1}(x)$}
\pscurve[linewidth=0.5pt,showpoints=false]{-}(5,.5)(5.3,1)(5.6,2.3)(6,3)
\pscurve[linewidth=0.5pt,showpoints=false]{-}(5.5,.6)(5.8,1.1)(6.1,2.4)(6.5,3.1)
\pscurve[linewidth=0.5pt,showpoints=false]{-}(6.8,.5)(7.1,1)(7.4,2.3)(7.3,3)
\pscurve[linewidth=0.5pt,showpoints=false]{-}(7.5,.8)(7.8,1.3)(8.1,2.6)(8.5,3.3)
\pscurve[linewidth=0.5pt,showpoints=false]{-}(5,.5)(5.5,.6)(6.8,.5)(7.5,.8)
\pscurve[linewidth=0.5pt,showpoints=false]{-}(5.3,1)(5.8,1.1)(7.1,1)(7.8,1.3)
\pscurve[linewidth=0.5pt,showpoints=false]{-}(5.6,2.3)(6.1,2.4)(7.4,2.3)(8.1,2.6)
\pscurve[linewidth=0.5pt,showpoints=false]{-}(6,3)(6.5,3.1)(7.3,3)(8.5,3.3)
\psline[linewidth=0.5pt]{->}(9.7,.5)(10,2.2)
\rput(9.5,1.4){$cs$}
\psline[linewidth=0.5pt]{->}(10,0)(11.5,0)
\rput(10.7,-.5){$cu$}
\pscurve[linewidth=0.5pt,showpoints=false]{-}(10,.5)(10.5,.6)(12.7,.5)(13.5,.8)
\pscurve[linewidth=0.5pt,showpoints=false]{-}(10,.7)(10.5,.8)(12.7,.7)(13.5,1)
\pscurve[linewidth=0.5pt,showpoints=false]{-}(10,.9)(10.5,1)(12.7,.9)(13.5,1.2)
\pscurve[linewidth=0.5pt,showpoints=false]{-}(10,1.1)(10.5,1.3)(12.7,1.1)(13.5,1.4)
\psline[linewidth=0.5pt]{-}(10,.5)(10,1.1)
\psline[linewidth=0.5pt]{-}(13.5,.8)(13.5,1.4)
\psline[linewidth=0.5pt]{-}(10.5,.6)(10.5,1.3)
\psline[linewidth=0.5pt]{-}(12.7,.5)(12.7,1.1)
\psline[linewidth=0.5pt]{->}(3.5,3)(4.5,3)
\rput(4,3.3){$g^n$}
\psline[linewidth=0.5pt]{->}(9,3)(10,3)
\rput(9.5,3.3){$g^n$}
\vspace{2cm}
\caption{}\label{fig5}
\end{flushleft}
\end{figure}
Finally, for condition (H3), a similar proof can be done as in the previous example.

\subsubsection{Example 4: Derived from Anosov through a Hopf bifurcation}

This example is the one treated in \cite{Mc1993}. Explicit formulas and details can be founded there. This example, obtained from a linear Anosov through a Hopf bifurcation is not generic and does not include  the examples in \cite{C1993}. A sharper analysis should be done to study the latter. 

Let $f:\T^3\to \T^3$ be a linear Anosov diffeomorphism, $T\T^3=E^s\oplus E^u$ with $\dim E^s=2$ and the $f$ has complex eigenvalues in $E^s(p),$ where  $p$ is  a fixed point. We will deform $f$ inside a small ball $B(p,r)$ to obtain a diffeomorphism $g:\T^3\to \T^3$ with the following features:
\begin{itemize}
\item $p$ is a repeller for $g.$ 
\item $g$ is partially hyperbolic $T\T^3=E^{cs}\oplus E^u.$
\item $Dg$ uniformly contracts $E^{cs}$ outside $B(p,r).$
\item $\|Dg_{/E^{cs}(x)}\|\le 1$ for any $x\notin W_{loc}^u(p).$
\item $\dist_{C^0}(g,f)<r.$
\item $E^{cs}$ is uniquely integrable.
\end{itemize}

\begin{figure}[h]
\begin{flushleft}
\vspace*{3.5cm}
\hspace*{1cm}
\psdot*[dotscale=.7](3,2)
\rput(3.5,1.5){$p$}
\psline[linewidth=0.5pt]{->>}(3,2)(1.5,1)
\pscurve[linewidth=0.5pt,showpoints=false]{->}(4,3)(3,2.7)(2.5,2)(3,1.6)(3.5,2)(3.4,2.2)(3,2.3)(2.8,2.1)
\psdot*[dotscale=.7](8,2)
\rput(8.2,1.5){$p$}
\psline[linewidth=0.5pt]{->>}(8,2)(6.5,1)
\pscircle[linewidth=0.5pt](8,2){.7}
\pscurve[linewidth=0.5pt,showpoints=false]{->}(8.2,3.5)(7.8,3)(7,2)(8,1)(9,2)(8,2.9)(7.8,2.8)
\pscurve[linewidth=0.5pt,showpoints=false]{->}(8,2.1)(7.85,2)(8,1.8)(8.2,2)(8,2.3)(7.6,2)
\psline[linewidth=0.5pt]{-}(2,1)(2,3.5)
\psline[linewidth=0.5pt]{-}(4.5,1)(4.5,3.5)
\psline[linewidth=0.5pt]{-}(2,1)(4.5,1)
\psline[linewidth=0.5pt]{-}(2,3.5)(4.5,3.5)
\rput(1,1.2){$E^u(p)$}
\rput(5.1,3.5){$E^s(p)$}
\rput(3.5,0.5){$A$}
\psline[linewidth=0.5pt]{-}(6.75,0.9)(6.75,3.5)
\psline[linewidth=0.5pt]{-}(9.25,.9)(9.25,3.5)
\psline[linewidth=0.5pt]{-}(6.75,.9)(9.25,.9)
\psline[linewidth=0.5pt]{-}(6.75,3.5)(9.25,3.5)
\rput(6,1.2){$E^u(p)$}
\rput(10.1,3.5){$E^{cs}(p)$}
\rput(8,0.5){$g$}
\caption{}\label{fig0}
\end{flushleft}
\end{figure}

We will prove that $g$ in the above conditions  has a unique probability measure of maximal entropy.
 Denote by $h$ the semiconjugation between $f$ and $g.$. It is not difficult to see that $h$ is injective on each $W^u(y,g)$ for any $y$ and moreover, $h(W^u(y,g))=W^u(h(y),f).$  Indeed, $h$ preserves the unstable foliation and the central stable foliation $h(W^{cs}(y,g))=W^s(h(y),f).$ Therefore,  an equivalent class is contained in a central stable manifold. 

\begin{lema}
Let $x\notin W^{cs}(p)$. Then $\diam (g^n[x]) \to 0.$
\end{lema}
\begin{proof}
Since $x\notin W^{cs}(p)$ then, there are infinitely many $n\ge 0$ such that  $g^n(x)\notin B(p,r)$ (and we may assume that without loss of generality that $g^n(x)\notin W^u_{loc}(p)).$ Therefore,
$\|Dg^n_{/E^{cs}(x)}\|\to 0$ and the same holds for any $y\in [x]$ (and uniformly on $y.$) The conclusion follows.
\end{proof}

\begin{cor}
Conditions {\rm (H1), (H2)} and {\rm (H3)} holds.
\end{cor}

\begin{proof}
First notice that the class $[p]$ is a closed disc, with $p$ as a repeller and the boundary attracts everything on the disk but $p.$ Therefore $h_{top}(g,[p])=0. $ If $[x]\subset W^{cs}(p)$ and $[x]\neq [p]$ then the class $[x]$ is attracted by the invariant circle and so $h_{top}(g,[x])=0. $ Now, if $[x]\notin W^{cs}(p)$, then $\diam (g^n[x])\to 0$ and therefore, for any $\epsilon$ and any $n$ large enough the cardinality of any $(n,\epsilon)$ separated set in $[x]$ is bounded, and hence $h_{top}(g,[x])=0. $ We have proved that (H1) holds.

For (H2), let $[x]$ be a periodic class. Notice that $[x]\cap W^{cs}(p)=\emptyset$ and therefore $\diam (g^n[x])\to 0.$ But since, $[x]$ is periodic, this means that $[x]=x$ and so $x$ is periodic and (H2) holds trivially.

Condition $(H3)$ can be proved with similar methods as the previous examples.
\end{proof}

\subsection{Suspension flows over derivated from Anosov diffeomorphisms}
Let $f\::\:X\to X$ be an homeomorphims defined on the compact metric space $X$. Let $\tau\::\:X\to\R^+$ be a continuous function and consider the space
$$Y=\{(x,t)\in X\times\R\::\:0\leq t\leq \tau(x)\},$$
with the points $(x,\tau(x))$ and $(f(x),0)$ identified for each $x\in X$. The {\em suspention flow} over $f$ with {\em height function} $\tau$ is the flow $\Phi=(\p_t)_{t\in\R}$ on $Y$ defined by
$$\p_t(x,s)=(x,s+t) \mbox{ whenever } s+t\in [0,\tau(x)].$$
Denote by $M(\Phi)$ and $M(f)$ respectively, the space of invariant measures of $\Phi$ and $f$. Then there exists a bijection $R\::\:M(f)\to M(\Phi)$ given by
$$R(\eta)=\frac{(\eta\times\m|Y)}{(\eta\times\m)(Y)}.$$
It is well known that
\begin{equation}\label{eq:entroflow}h_{top}(\Phi)=\sup_{\xi\in M(\Phi)}h_\xi(\p_1),\end{equation}
where $\p_1$ is the time one map of the flow $\Phi$.
Setting $R(\eta)=\xi$, by the Abramov's formula \cite{Ab1959}, we have
\begin{equation}\label{eq:Abramov}
 h_{\xi}(\p_1)=\frac{h_{\eta}(g)}{\eta(\tau)}.
\end{equation}

Bowen and Ruelle \cite{BoRu1975} showed that in the particular case when $f$ is a  transitive Anosov diffeomorphism and the height function $\tau$ is H\"older, then there exist a unique measure $\mu_\Phi$ maximizing the entropy of $\Phi$ and
$$\mu_\Phi=R(\mu_\p),$$
where $\mu_\p$ is the unque equilibrium state of the potential $\p=-h_{top}(\Phi)\tau$. In the particular case, when the height function $\tau$  is cohomological to a constant, then $\mu_\p=\mu$ the Bowen measure for $f$.
Now consider $g\::\:X\to X$ and homeomorphism satisfying the conditions in Section~\ref{sec:statements} and let $\Psi$ be the suspension flow over $g$ with the same height function $\tau$, H\"older continuous and cohomological to a constant. By simplicity, we will take $\tau(x)=1$ for all $x\in X$. The general case is studied in \cite{Samva2009}.
\begin{cor}\label{cor:ultimo}
The flow $\Psi$ has a unique measure maximizing the entropy $\mu_\Psi=\mu_g\times \m$, where $\mu_g$ is the Bowen measure for $g$.
\end{cor}
\proof First, note that $R(\mu_g\times \m)=\mu_g$ so from \eqref{eq:Abramov} and since $\mu_g(\tau)=1$ we have
\begin{equation}\label{eq:entropflu}h_{\mu_g\times \m}(\psi_1)=h_{\mu_g}(g)=h_{top}(g).\end{equation}
Considering \eqref{eq:entroflow}, \eqref{eq:Abramov} and the bijection betwen $M(\Psi)$ and $M(g)$, we have
\begin{equation}\label{eq:entropflu2}h_{top}(\Psi)=\sup_{\xi\in M(\Psi)}h_\xi(\p_1)=\sup_{\eta\in M(g)}h_{\eta}(g)=h_{top}(g).\end{equation}
Finally from \eqref{eq:entropflu} and \eqref{eq:entropflu2} we obtain that $\mu_g\times \m$ is the unique measure maximizing the entropy of $\Psi$.

\endproof

Of course, Corollary~\ref{cor:ultimo} implies that all the suspension flows over each derivated from Anosov diffeomorphisms considered in the examples above, with height function a constant have unique measure maximizing the entropy.

\bibliographystyle{plain}
\bibliography{SAMVA.bib}
\end{document}